\documentclass{siamart220329}
\usepackage{lipsum}
\usepackage{amsfonts}
\usepackage{graphicx}
\usepackage{epstopdf}
\usepackage{algpseudocode}
\usepackage{caption}
\usepackage{subcaption}
\usepackage{tikz-cd}
\usepackage{pgfplots}
\usepackage{wrapfig}
\pgfplotsset{compat=1.18}
\usepackage{titlesec}

\usepackage{caption}
\usepackage{enumitem}
\setlist{noitemsep,topsep=1pt}

\ifpdf%
  \DeclareGraphicsExtensions{.eps,.pdf,.png,.jpg}
\else
  \DeclareGraphicsExtensions{.eps}
\fi
\usepackage{amsopn}

\usepackage{booktabs}

\def\Popt{P_{\sharp}}      
  
\def\Pid{P_{\star}}      

\def\XNorm#1{\left\| #1 \right\|}                     
                                  
\def\XRange#1{{\rm range}(#1)}  
                                
\def\XVec#1{{\mathbf #1}}       

\def\Xu{\XVec{u}}
\def\Xv{\XVec{v}}

\def\Xf{\XVec{f}}
\def\Xe{\XVec{e}}
\def\Xr{\XVec{r}}
\def\Xn{\XVec{n}}

\def\Xdof{\mathbf{dof}}

\def\XfG{\mathbf{f_N}}
\def\XcG{\mathbf{c_N}}

\newcommand{\ZQS}{{Q_S}}
\newcommand{\ZR}{{R}}
\newcommand{\ZQR}{{Q_R}}
\newcommand{\ZS}{{S}}

\def\XReals{{\mathbb R}}

\newcommand{\XMt}{{\widetilde{M}}}

\def\XMt{{\widetilde{M}}}
\newcommand{\Hdiv}{{H(\mbox{div})} }            
\newcommand{\Hcurl}{{H(\mbox{curl})} }

\newcommand{\TheTitle}{%
 Nodal AMG Coarsening and Interpolation for PDE systems
}

\newcommand{\TheFunding}{%
  We thank Randy Bank for his valuable insight and guidance, without which these results would not fully come to life.
}

\author{James Brannick\thanks{Department of Mathematics, The Penn State University, University Park, PA 16802, USA (\email{jjb23@psu.edu})} \and Robert Falgout\thanks{Lawrence Livermore National Laboratory, Livermore, CA, 94550, USA (\email{falgout2@llnl.gov})} \and Karsten Kahl\thanks{Fachgruppe Mathematik \& Informatik, Bergische Universität Wuppertal,  Wuppertal, North Rhine-Westphalia, 42119, Germany (\email{kkahl@uni-wuppertal.de})}\and Jacob Schroder\thanks{Department of Mathematics \& Statistics, The University of New Mexico, Albuquerque, New Mexico, 87106, USA (\email{jbschroder@unm.edu})}\and Taoli Shen\thanks{Department of Mathematics, Penn State University, University Park, PA 16802, USA (\email{tzs5765@psu.edu})}}
\title{{\TheTitle}}
\title{{\TheTitle}\thanks{\TheFunding}}

\begin{document}

\maketitle


\begin{abstract}
  We present an approach to constructing a practical coarsening algorithm and interpolation operator for the algebraic multigrid (AMG) method, tailored towards systems of partial differential equations (PDEs) with large near-kernels, such as $\Hcurl$ and $\Hdiv$. Our method builds on compatible relaxation (CR) and the ideal interpolation model within the generalized AMG (GAMG) framework but introduces several modifications to define an AMG method for PDE systems. We construct an interpolation operator through a coarsening process that first coarsens a nodal dual problem and then builds the coarse and fine variables using a matching algorithm.  Our interpolation follows the ideal formulation; however, we enhance the sparsity of ideal interpolation by decoupling the fine and coarse variables completely. When the coarse variables align with the geometric refinement, our method reproduces re-discretization on unstructured meshes. Together with an automatic smoother construction scheme that identifies the local near kernels, our approach forms a complete two-grid method. Finally, we also show numerical results that demonstrate the effectiveness of this interpolation scheme by applying it to targeted problems and the Stokes system. 

  
\end{abstract}


\section{Introduction}\label{sec:intro}
AMG is an algorithm for solving linear systems of equations typically derived from the discretization of a partial differential equation,
\begin{equation}
\label{eqn:axb}
A \Xu = \Xf,
\end{equation}
where $A \in \mathbb{R}^{n\times n}$ is a sparse matrix. AMG constructs a fast solver through a complementary process of relaxation and coarse-grid correction. Unlike geometric multigrid (GMG), AMG and adaptive and bootstrap AMG methods are designed to treat smooth as well as oscillatory errors that yield small residuals after relaxation. The basic approach in AMG is to identify the smooth error (i.e., error not effectively reduced by relaxation) and then coarsen and build interpolation that directly approximates and treats such error.  For scalar systems, one assumes the constant represents smooth errors locally, and this assumption is the key to the design of classical AMG~\cite{RuStu1987}. 

For PDE systems, the notion of smooth error and complementarity of coarse-grid correction (for a fixed smoother) is not as well understood.
Our focus is on developing AMG techniques and obtaining insights into the AMG coarsening process for PDE systems.
The two-grid scheme is summarized by the error propagator
\begin{equation}\label{eq:tl}
E_{TG} = (I-M^{-1}A)(I - P A_c^{-1} P^T A),
\end{equation} 
where $M$ defines the smoother, $P\in \mathbb{R}^{n \times n_c}$ denotes the interpolation matrix that maps error corrections to the fine level and $A_c = P^TAP$ is the coarse-level system matrix. 
One of the main goals of this paper is to define a proper splitting of the degrees of freedom (DoFs) into coarse and fine.

Assuming that the system matrix $A$ is symmetric, the generalized AMG two-grid theory can be summarized as follows~\cite{FalgoutGAMG,sharp_theory_2005,G-BAMG-2018}.
\begin{equation} \label{eqn-tg-measure}
\begin{split}
  \XNorm{E_{TG}}_A^2 \leq 1 - \frac{1}{K} , &~~~\mbox{where}~~~
  K := K(PR) = \sup_{\Xe} \frac{ \XNorm{(I - PR)\Xe}_{\XMt}^2 }{ \XNorm{\Xe}_A^2 } \geq 1 ,
\end{split}
\end{equation}
where $\XMt = M(M+M^T-A)^{-1}M^T$ is the symmetrized smoother. Here, $R: \XReals^{n} \mapsto \XReals^{n_c}$ is any matrix for which $RP = I_c$, the identity on $\XReals^{n_c}$, so that $PR$ is a projection onto $\XRange{P}$. 
The derivation of the GAMG theory~\cite{FalgoutGAMG} which we use to design our AMG solver 
 begins with the $\ell_2$-space decomposition of $V=\mathbb{R}^n =\mathbb{R}^{n_s} \oplus\mathbb{R}^{n_c}$ into a basis for the fine variables, the columns of the matrix $S\in\mathbb{R}^{n\times n_s}$, and a basis for the coarse variables, the columns of the matrix $R^T \in \mathbb{R}^{n\times n_c}$. If we assume the orthogonality conditions  $RS=0$, $S^TS = I_s$, and $RR^T = I_c$ as in GAMG, then 
\begin{equation}\label{eqn:split}
    \Xu = \ZQS \Xu + \ZQR \Xu = S\Xu_s + R^T \Xu_c, \quad \ZQS = \ZS \ZS^T, \quad \ZQR = \ZR^T\ZR, 
\end{equation}
with $\Xu_s$ denoting the fine variables and $\Xu_c$ the coarse variables in these new bases.  The corresponding generalized form of ideal interpolation used as  motivation is then:
\begin{equation} \label{eqn-ideal}
\begin{split}
  P_{\star} := \arg\min_{P:\, RP=I_c} K(PR)
  &=
  \left[~ S ~~ R^T ~\right]
  \left[ \begin{array}{c}
    - (S^T A S)^{-1} (S^T A R^T) \\ I
  \end{array} \right] 
  \\[1ex]
  &=
  (I - S(S^T A S)^{-1} S^T A) R^T.
\end{split}
\end{equation}
This form of ideal interpolation gives a global harmonic extension in these new bases that minimizes the convergence rate for a fixed smoother. Importantly, the ideal interpolation operator is closely related to the convergence of CR, namely, fast convergence of CR implies uniform convergence of the two-level method using ideal interpolation.  

A general two-level method in the GAMG framework using CR is given by a smoother $B\approx A^{-1}$ and choices of $R$ and $S$ that define the coarse-level correction, yielding the error propagator:
\begin{equation}\label{eq:GAMGtl}
    E_{TL} = (I - \ZQS \: B \: \ZQS A)(I-\pi_A(\Pid)), \quad \pi_A(\Pid) := \Pid (\Pid^T A \Pid)^{-1} \Pid A,
\end{equation}
where the first term defines the CR iteration.  
 Assuming that $A$ is SPD and the smoother is invertible, we define the primary compatible relaxation iteration and its error transfer matrix, both based on this $\ell_2$ splitting, as follows:
 \begin{equation}\label{CR}
\Xu^{k+1} = \Xu^k  + \ZQS M^{-1} \ZQS A \Xr_k
\quad \mbox{and} \quad 
\Xe^{k+1} = ( I - \ZQS M^{-1}\ZQS A)\Xe^k,
\end{equation}
where we set $B=M^{-1}$.
We note that other CR methods are possible as error transfer iterations; for example, we can define $S$-relaxation and Habituated CR, respectively:
\begin{equation}\label{CR-error2}
\Xe_\ZS^{k+1} =  (I_s - M_\ZS^{-1} A_\ZS)\Xe_\ZS^k \quad \mbox{and} \quad 
\Xe^{k+1} = \ZQS( I - M^{-1} A)\ZQS\Xe^k,
\end{equation}
where $X_\ZS := \ZS^T X \ZS$ for an invertible matrix $X$.  The former $S$-relaxation method can be turned into an actual iterative method with proper choices of $\ZR$ and $\ZS$ as shown below, whereas the latter habituated CR accounts for the global smoother's action and, hence, can be used to tune the smoothing parameter to obtain fast convergence.

\section{Featured systems}
Our focus in this paper is to develop AMG methods for near-singular PDE systems of the bilinear form:
\begin{equation} \label{eqn-bilinear}
    A(\Xu, \Xv) = (\Xu,\Xv)_X+\beta(\Xu, \Xv),
\end{equation}
for vector function $\Xu, \Xv \in X$, where $ \left( \cdot, \cdot \right),  \left( \cdot, \cdot \right)_X$ are the inner product on $L^2(\Omega)$ and $X$, $X = H(\Omega;\mbox{curl})$ or $H(\Omega;\mbox{div})$ and $\beta >0$ is small.
When $\beta \xrightarrow{} 0$, these problems have large near-kernels which make traditional pointwise smoothers ineffective. 
For example, in the case of $\Hcurl$, equation \eqref{eqn-bilinear} becomes the shifted curl-curl system. 
Discretizing it with lowest-order N\'{e}d\'{e}lec edge finite elements~\cite{nedelec1980mixed, hiptmair2002finite} gives
$A = A_s + \beta A_m \in \XReals^{n \times n}$ with $A_s$ and $A_m$ denoting the stiffness and the mass matrices. An important property of the curl-curl system is that the resulting matrix $A$ has a large $O(n)$ local near nullspace when $\beta$ is small.
In this case, non-pointwise smoothers~\cite{hiptmair1998multigrid} are needed in order to obtain an optimal-order multigrid method.  

\section{Automated construction of smoothers}
The first task in designing AMG for PDE systems is to define an appropriate smoother.  For problems with local near kernels, we propose the use of the algorithm in~\cite{claus2019multigrid} for the automated construction of AMG smoothers.  The method simultaneously builds the smoother and computes the local near kernel that we use to define the coarse variables.  Given its central role in our AMG algorithm, we include a detailed description of the method from~\cite{claus2019multigrid}.   

An automatic scheme for smoother construction is proposed in~\cite{claus2019multigrid}. Given $A$ and a coarsening factor $m$ (used in the AMG coarsening algorithm to select the coarse variables), the algorithm starts by constructing all {\em diameter} $m$ sets $D_i$ for each DoF $i$. Next, we check for near null space components in the corresponding principal submatrix by computing a local eigen-decomposition (or more generally an SVD). To construct a diameter $m$ set for a given DoF $i$, we first find all its distance $m$ neighbors based on the matrix graph of $A$. For each j in this neighborhood, we define the diameter $m$ set as the intersection of all distance $k$ neighborhoods of $i$ and $j$, where $k < m$. The overall approach is summarized in algorithm~\ref{alg-G}. The parameter $\epsilon$ serves as a filter to select only the near-kernel modes.

\begin{algorithm}
\caption{Computing local near kernel of $A$ and associated index sets}
    \begin{algorithmic}[1]
    \algnewcommand{\Initialize}[1]{%
    \Statex \hspace*{\algorithmicindent}\parbox[t]{.8\linewidth}{\raggedright #1}
    }
    \Function {findLocalNearKernels}{A, m, $\epsilon$}
    \Initialize {$N = \emptyset$}
    \For{all $i$}
    \State Construct list $D_i$ of all diameter $m$ sets for node $i$
       \For{all $\Omega_{ij}$ in $D_i$}
        \State Solve $A(\Omega_{ij},\Omega_{ij}) \Xv = \lambda_{\min} \Xv $
        \If{$\lambda_{\min}  \leq \epsilon$}
        \State $N = \left[
        \begin{array}{ccc}
        N, \, \Xu
        \end{array}
        \right], \text{ where }\Xu_{\Omega_{ij}} = \Xv \text{ and }0 \text{ elsewhere} $
        \EndIf
        \EndFor
    \EndFor
    \State \Return $N$
    \EndFunction
    \end{algorithmic}
    \label{alg-G}
\end{algorithm}
Given the near-kernel matrix $N$ and its supports defined by index sets $\{\Omega_\ell\}_{l=1}^L$, we consider distributive and overlapping multiplicative Schwarz methods
\begin{equation} \label{eqn-g-smoother}
 E_D =  I - N (N^T A N )^{\sim 1} N^T A \quad \mbox{and} \quad E_{MS} = \prod_{i=1}^{L} \left( I - I_i A_i^{-1} I_i^T A \right ),
\end{equation}
respectively. The notation $\sim\!1$ for the distributive smoother indicates ``approximate inverse'' using one sweep of a pointwise smoother, e.g. Gauss-Seidel, to ensure global smoothing on all 
DoFs.  
Overall, the algorithm leads to distributive smoothers such as the Hiptmair smoother or star-patch relaxation methods based on the Pavarino-Arnold-Faulk-Winther type splittings in the case $N$ corresponds to discrete grads, curls, and divs, respectively.

\section{Nodal coarsening and building AMG interpolation}\label{sec:interpolation}
We proceed with a description of the two-level setup algorithm, i.e., our coarsening algorithm. 
 The following three steps summarize the overall approach:
 \begin{enumerate}
\item Use the Algorithm~\ref{alg-G} to compute the local near kernel matrix $N$.   
 \item 
 Coarsen $A_N = N^TAN$ using classical C/F splitting of $\XcG$ and $\XfG$, which correspond to indices in the local near kernels. Then build appropriate coarse and fine variables $R^{T}$ and $S$ by Algorithm~\ref{alg:p}.
 \item 
 Compute the operator defined in~\eqref{eqn-ideal} and the Galerkin operator $P^TAP$.
\end{enumerate}
In the case of $H^1$, we set $N$ as the constant vector and run classical AMG. In practice, we can use CR to monitor the quality of the coarse grid as fast convergence of CR implies uniform convergence of the two-level method. Hence, CR not only guides our coarsening process but also serves as a main theoretical tool in the design of AMG.  
 
\subsection{AMG interpolation}
The construction of our AMG interpolation uses the definition of the generalized ideal interpolation operator defined in~\eqref{eqn-ideal}, 
where, again, $R$ and $S$ store the coarse and fine variables as in~\eqref{eqn:split}.
In classical AMG, the coarse and fine variables take the form
\begin{equation}\label{eq:classicalRS}
 \quad R= \left[\begin{array}{cc}  0  & I \end{array}\right] \\ 
, \quad S = \left[\begin{array}{c} I \\
 0 \end{array}\right]
\end{equation}
It is clear that this type of coarsening simply injects a subset of fine DoFs into the coarse mesh without preserving any near kernels structures. This becomes particularly problematic when local near kernels overlap and errors consisting of the linear combination of these overlapping near kernel components lead to stagnation in the convergence of the two-grid method with standard ideal interpolation.  Instead, for these problems, it is important to preserve the proper averaging (orientation) of the local near kernel components when defining the coarse variables.

\subsection{Inspiration from geometric interpolation}
\begin{figure}[t]
    \centering
    \begin{subfigure}[t]{0.24\textwidth}
        \centering
        \includegraphics[width=\linewidth]{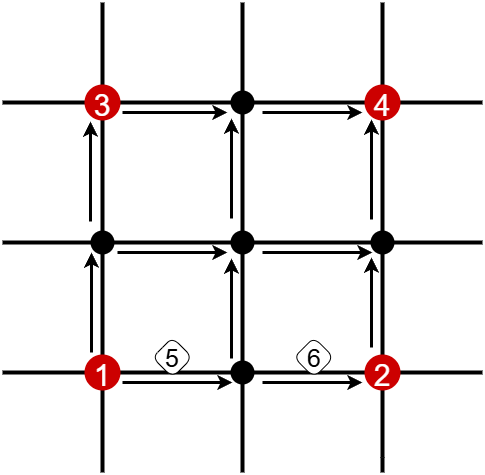}
        \caption{Nodal coarsening}
        \label{fig:R_orientation}
    \end{subfigure}
    \hfill
    \begin{subfigure}[t]{0.24\textwidth}
        \centering
        \includegraphics[width=\linewidth]{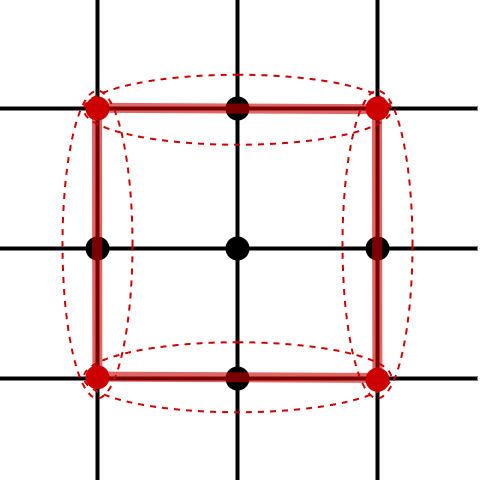}
        \caption{Pattern in $R$}
        \label{fig:R_quad}
    \end{subfigure}
    \hfill
    \begin{subfigure}[t]{0.24\textwidth} 
        \centering
        \includegraphics[width=\linewidth]{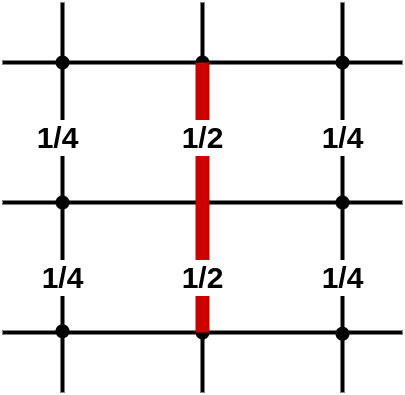}
        \caption{Vertical stencil}
        \label{fig:vertical_stencil}
    \end{subfigure}
    \hfill
    \begin{subfigure}[t]{0.24\textwidth} 
        \centering
        \includegraphics[width=\linewidth]{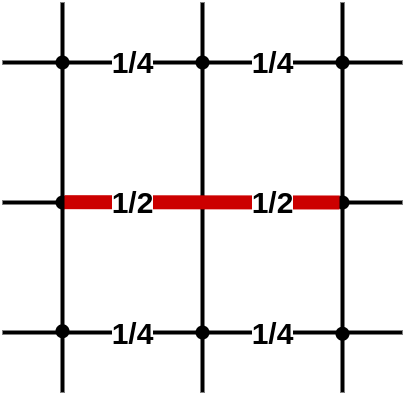}
        \caption{Horizontal stencil}
        \label{fig:horizontal_stencil}
    \end{subfigure}
    \caption{Curl-curl on uniform quadrilateral mesh.}
    \label{fig:curl-curl}
\end{figure}

To inspire the algebraic construction of our interpolation, we begin with the curl-curl problem on a uniform quadrilateral mesh. A window of the mesh of 4 elements with 12 DoFs, each defined by the tangential component of an edge, is illustrated in Figure~\ref{fig:R_quad}. We can partition these 12 DoFs into the external DoFs highlighted in red and the interior DoFs. The corresponding geometric coarse mesh, also highlighted in red, contains 4 DoFs, represented by the dashed red circles enclosing the fine DoFs. Specifically, each coarse edge is subdivided into two fine edges, and the mesh spacing between the nodes is doubled. Here, the red circles also depict our matching coarsening algorithm.  

The overall coarsening we choose is then described by the choice of coarse and fine variables on the exterior and interior DoFs. To explain our findings, we begin with the geometric choice of $R$ and $S$:
\begin{equation}\label{eq:RS}
  R_E := \begin{array}{c}
  \frac{1}{\sqrt{2}}\left[\begin{array}{ccc}  I & I \end{array}\right]
  \end{array}, \quad 
  S_E := \frac{1}{\sqrt{2}}\left[\begin{array}{c} 
  I  \\
  -I 
  \end{array}\right] 
  ,  \quad \mathrm{and} \quad S_I := \left[\begin{array}{c} I \\
 0 \end{array}\right] ,
\end{equation}
with subscripts $I,E$ denoting the interior and exterior edges, respectively. Here, we notice that the choice of $R_E$ and $S_E$ coincides with that of Theorem 6.2 in ~\cite{FalgoutGAMG}. Although the example there is $H(\operatorname{div})$, the geometric coarse mesh is identical. 

In our construction, we take the AMG approach and instead build the near kernel that is encoded in $N$ directly into the definition of the coarse variables so that 
\begin{equation}\label{eq:newRgrad}
  R_{E_G} := \begin{array}{c}
  \frac{1}{\sqrt{2}}\left[\begin{array}{ccc}  \pm I & \pm I \end{array}\right]
  \end{array}.
\end{equation}
Here, the $+$ or $-$ must be consistent with the underlying orientations in the fine grid near kernels. For instance, in Figure~\ref{fig:R_orientation}, since edges 5 and 6 share the same orientation, the corresponding row in $R$ contains 1's at positions 5 and 6, and 0's elsewhere. If one of the edges has an opposite orientation, its sign in $R$ must be flipped accordingly. The coarse variables now preserve the near kernels (grad, curl and/or div) exactly. By following correctly oriented paths, we can algebraically construct arbitrary coarse near-kernel components. 

Next, we adopt the AMG approach by incorporating the near kernels directly into the coarse variables. This is achieved by moving $S_E$ to the coarse variables:
\begin{equation}\label{eq:newRgradglobal}
  R := \left[\begin{array}{ccc}  R_{E_G}  \\ S_E^T \end{array}\right] \quad \text{and} \quad S := S_I.
\end{equation}
Numerically, we observe that the new coarse variables generate a coarse space that preserves all necessary gradients and curls. For CR, $S^TAS$ now is just $A_{FF}$, and hence, the AMG smoother we use is fixed. Therefore, our $S-$relaxation CR defined in~\eqref{CR} is a well-defined reduction-based smoother and iteration~\eqref{eq:GAMGtl} defines a general two-level reduction-based AMG method, where $F-C$ smoother can be used in practice.

Finally, we consider computing a sparse approximation to our new ideal $P$.  Notably, we find that we can completely drop the $S_E$ term from the definition of $R$ and that the range of this approximation and that of using $R_E$ and $S_E$ defined in \eqref{eq:RS} are the same.  Though we have not found an appropriate analytical form, we are able to compute a post-scaling transformation numerically that establishes the equivalence between the range of these two interpolations for a set of example problems and coarse grids. One explanation that we do not need $S_E$ in our coarse variables is that the smoother has already dealt with it. The final choices of $R$ and $S$ that we use in our algorithm are now given by  
\begin{equation}\label{eq:newRgradglobal2}
  R :=  R_{E_G} \quad \mbox{and} \quad S := S_I,
\end{equation}
which gives the highly localized choices that at the same time preserve the near kernels. 
Given $R$, we compute local harmonic extensions of the tentative interpolation operator  $R^T$, extending the averages to the DoFs defined by $S$. We note that, by following geometric mesh refinement to define the coarse variable in $R^T$, we are able to reproduce re-discretization and the boundary conditions on the coarse grid for structured as well as unstructured meshes.
Our new splitting does not follow the GAMG theory in that our definitions of $R$ and $S$ in the approximation we construct only provide a partial splitting of $\mathbb{R}^n$.  Specifically, we have 
$\dim(\text{span}(R^T)) + \dim(\text{span}(S)) = n_c + n_s \leq n.
$
We note that, our derivation using AMG techniques recovers the theoretical result in~\cite{Marian} where the dependence on the partition of unity is removed for AMGe interpolation. 

\subsection{Algebraic construction of our coarse and fine variables} \label{sec:algebraic_p}
So far, in the design of our AMG interpolation operator, we have followed the mesh refinement. To make this process algebraic, we rely on the coarsening of the nodal dual space defined by $N^TAN$, assuming that $N$ is readily available from Algorithm~\ref{alg-G}. Conceptually, $N^TAN$ defines the scale of the local near kernels, determined by the diameter $m$ in Algorithm~\ref{alg-G}. Subsequently, we split the DoFs into exterior and interior groups. The complete algorithm is summarized below.
\begin{algorithm}[H] 
\caption{construct the interpolation $P$}\label{alg:p}
\begin{algorithmic}[1]
\algnewcommand{\Initialize}[1]{%
\Statex \hspace*{\algorithmicindent}\parbox[t]{.8\linewidth}{\raggedright #1}
}
\Function{form$P$}{$A$, $N$}
\Initialize{$S = \emptyset$, $R = \emptyset$, $\Xdof=\{1, \dots, n\}$}
\State{$\XcG, \XfG \gets \text{C/F splitting}\,(N^TAN)$} 
\For{$i, j \in \XcG$} \label{alg:split-start}
    \If{$i$ and $j$ are distance-two neighbors}
        \State $k \gets \text{mutual distance-one neighbor of } i \text{ and } j$
        \State Select $\{ DoF_1, DoF_2 \} \subset \Xdof$ that connect $i$ and $j \text{ through } k$
        \State $\Xdof \gets \Xdof \setminus \{ DoF_1, DoF_2 \}$
    \EndIf
    \State{$R^T =  \begin{bmatrix}
            R^T & \begin{bmatrix}
             \cdots& 0 &[\pm 1]_{DoF_1} &\cdots & 0 &\cdots  &[\pm1]_{DoF_2} & 0 & \cdots
        \end{bmatrix}^T
        \end{bmatrix}$}
\EndFor \label{alg:split-end}
 \State {$S \gets \begin{bmatrix}
            S & \begin{bmatrix}
            I_\Xdof & 0
        \end{bmatrix}^T
        \end{bmatrix}$} 
\State \Return {$P_\star = (I - S(S^T A S)^{-1} S^T A) R^T$}
\EndFunction
\end{algorithmic}
\end{algorithm}



In particular, line~\ref{alg:split-start} - line~\ref{alg:split-end} constructs the coarse variables in $R^T$. For every pair of coarse nodes $i,j$ distance-two neighbors, we connect them by averaging the two DoFs $DoF_1$ and $DoF_2$ that form a path between them. If multiple paths exist between them, we only pick one of them based on the labeling order for simplicity.  The definition of coarse variables must be consistent with the near kernel at the related fine DoFs, and the orientations of the near kernel for the fine DoFs can be determined from the corresponding rows of $N$.  The paths in the coarse variables naturally define the exterior DoFs. Therefore, the remaining DoFs belong to the interior where we assign them as Euclidean basis in $S$. In terms of building AMG interpolation, our approach shares similar ideas with~\cite{Reitzinger&Schöberl} to preserve near kernels on the coarse mesh algebraically through averaging with the correct orientations.

To illustrate the idea, we apply Algorithm~\ref{alg:p} to the same curl-curl example in Figure~\ref{fig:R_orientation}. First, we run a C/F splitting on $N^TAN$, which returns a nodal full coarsening set marked in red. We then form length-two paths by selecting edge DoFs between node pairs $\{1,2\},\{1,3\}, \{2,4\}, \{3,4\}$. As a result, coarse variables in $R^T$ are identical to those obtained by following the refinement pattern illustrated in Figure~\ref{fig:R_quad}. 

\begin{figure}[h]
    \centering
    \begin{subfigure}[t]{0.22\textwidth}
        \centering
        \includegraphics[width=\linewidth]{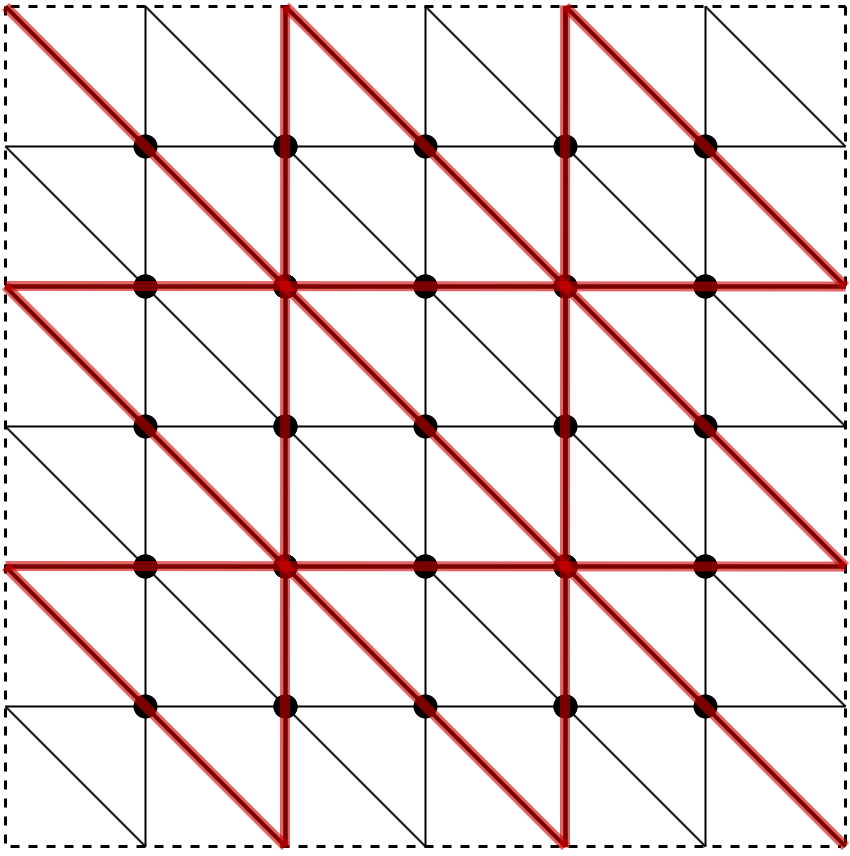}
        \caption{Refinement}
        \label{fig:triangle_refinement}
    \end{subfigure}
    \hfill
    \begin{subfigure}[t]{0.22\textwidth}
        \centering
        \includegraphics[width=\linewidth]{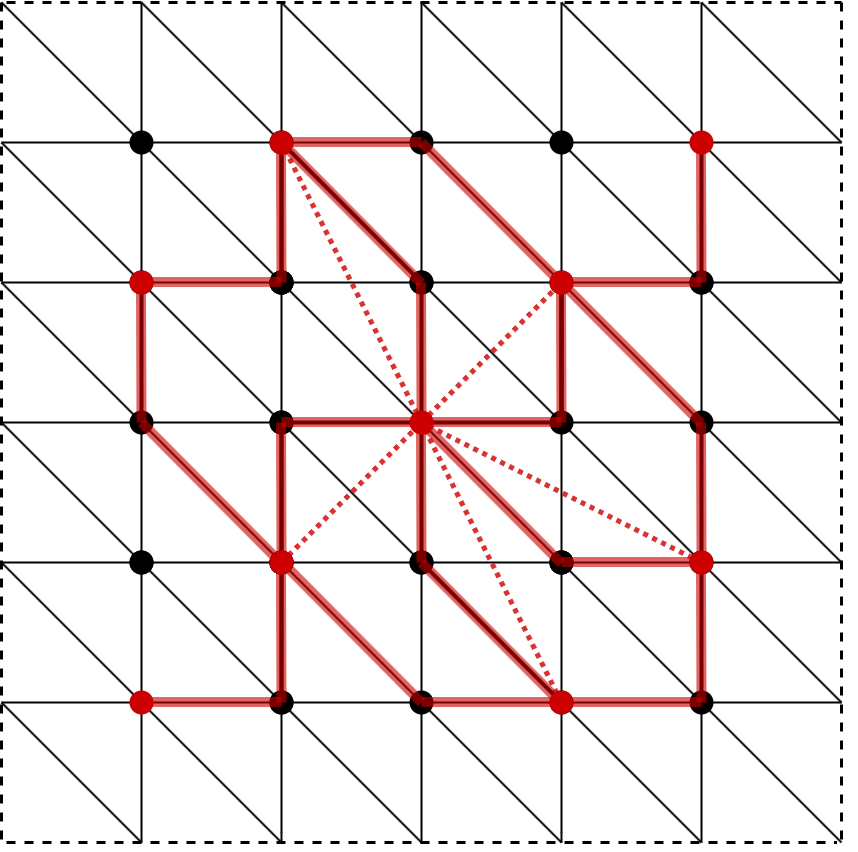}
        \caption{Algebraic pattern I}
        \label{fig:triangle_algebraic}
    \end{subfigure}
    \hfill
    \begin{subfigure}[t]{0.22\textwidth} 
        \centering
        \includegraphics[width=\linewidth]{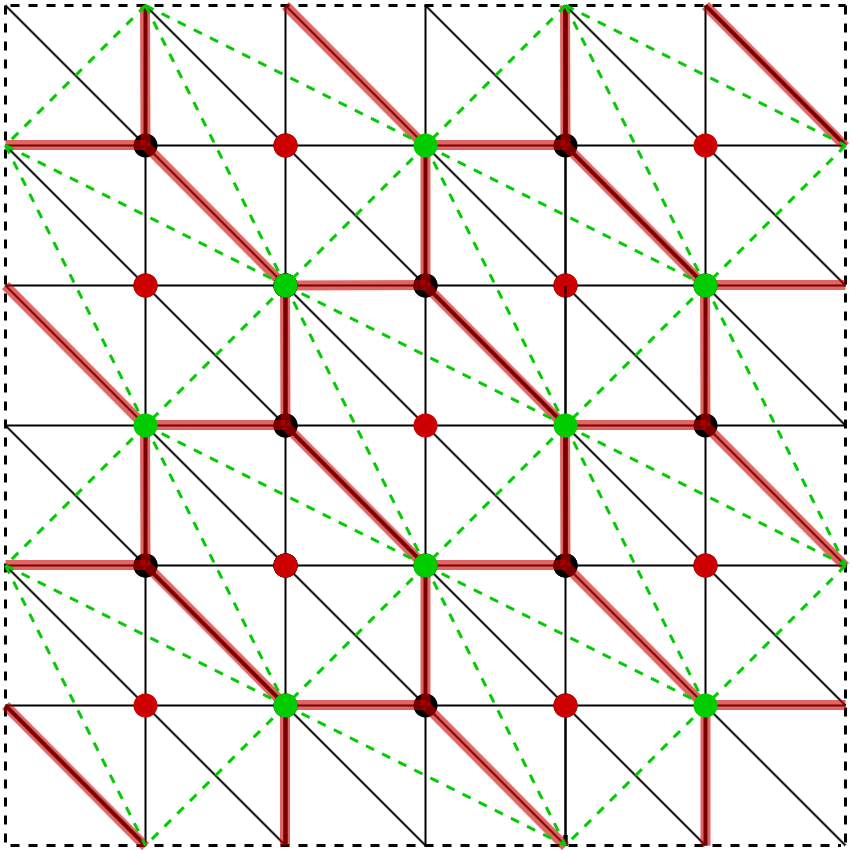}
        \caption{Algebraic pattern II}
        \label{fig:triangle_alternative}
    \end{subfigure}
    \caption{Uniform triangular mesh with Dirichlet boundary conditions}
    \label{fig:curlcultriangle}
\end{figure}

In general, Algorithm~\ref{alg:p} does not obtain geometric refinement. Nevertheless, we still maintain a near kernel structure algebraically on the coarse mesh. For example, consider the curl-curl problem on a uniform triangular mesh in Figure~\ref{fig:curlcultriangle}. Here, the Dirichlet boundary condition ($\Xu \times \Xn = 0$) is enforced along the black dotted line. Applying Algorithm~\ref{alg:p} produces the pattern shown in Figure~\ref{fig:triangle_algebraic}, where the length-two paths are highlighted as solid red straight lines. Notably, the coarse DoFs of length-two paths form a variant of $2h$-near kernels on the coarse mesh. One such $2h$-kernel is marked with red dotted lines, where each line represents a path between two coarse nodes. Due to the Dirichlet boundary condition, the absence of coarse nodes on the boundary prevents forming paths to the boundary. As a result, some near kernels near the boundary may be excluded when constructing coarse variables.

Figure~\ref{fig:triangle_alternative} presents an alternative coarsening strategy that can be executed algebraically. The idea is to take the support of the near kernels obtained from Algorithm~\ref{alg-G} at $\XcG$ as the interior, while the remaining DoFs, highlighted by red lines, are assigned to the exterior. The underlying principle remains the same. We form length-two paths between coarse nodes; however, the paths are formed between the green nodes, while the red coarse nodes are reclassified as fine nodes. This approach allows us to construct a uniform coarse mesh, represented by green dotted lines, that closely resembles geometric refinement. While this splitting of the interior and exterior can be obtained algebraically, a key challenge lies in the automatic construction of $R$ due to the overlapping paths in a specific pattern. A possible solution is to define a multi-vector preserving set of coarse variables that overlap. Overall, this approach offers another promising potential strategy for the design of a practical multilevel recursion that we intend to investigate.

\section{Solving saddle-point systems using curl-curl solvers: the Stokes problem}
\label{sec:stokes}
In solving saddle-point systems using our proposed AMG solver for PDE systems, we follow the same approach as before, except that now we assume the block structure of the PDE is given.  In this case, our algorithm closely resembles the methods developed in~\cite{Tong} for curl-curl problems, where the use of classical AMG coarsening for the nodes and aggregation coarsening for the edges is proposed.  Hence, we apply our setup algorithm and compute a block-diagonal interpolation operator $P$ for 2x2 block systems of the form:
   \begin{equation}\label{eq:A saddle}
  A := \left[\begin{array}{ccc}  A_{e}  & A_{e,N}  \\ A_{N,e} & 0 \end{array}\right]\quad , \quad P = \left[\begin{array}{ccc}  P_{e}  &   \\ & P_N\end{array}\right],
 \end{equation}
 with subscripts $e,N$ denoting edge and nodal DoFs as in curl-curl.  For the Marker and Cell (MAC) finite difference approximation to the Stokes system, we use the same approach that we used for the curl-curl problem. To justify this, consider the problem defined on a dual mesh with vertices located at the pressure (nodes) locations and velocities (edges) connecting them.  Since the element-DoF relations for the MAC scheme are identical to the element-DoF relations for the curl-curl system given in Figure~\ref{fig:curl-curl}, preserving velocities is again accomplished by averaging, i.e., matching edges between coarse nodes.  Hence, the coarse variables for the Stokes problem are the same as the ones we use to coarsen curl-curl.  To further justify this concept, we focus only the primary global Stokes system by assuming periodic boundary conditions, as in local Fourier analysis. Then, we apply the classical AMG coarsening algorithm to $$A_N = A_{N,e} A_e A_{e,N} = G^T A_e G,$$ where for Stokes we assume that $A_{N,e} = G^T,$ $A_{e,N} = G$, where $G$ is the discrete gradient, and we used that 
  $A_e = -\Delta + \beta I = \nabla \times \nabla \times - \nabla \nabla \cdot + \beta I$ with the assumption that $\nabla \cdot u =0$ and $\beta = 0$. In this view, up to boundary conditions,  coarsening the Stokes problem is equivalent to coarsening the \Hcurl problem from an AMG perspective, an observation recently made by Hiptmair for no-slip boundary conditions and Scott for Freudenthal meshes in 2d and 3D~\cite{boon2024hcurlbasedapproximationstokesproblem,doi:10.1137/22M1533943}. 
Finally, we mention that, if we use lowest-order Raviart-Thomas FEM for grad-div systems, then our solver works for this problem without any modification.  Hence, overall, in our AMG framework, the curl-curl problem, the grad-div problem, and Stokes' system are essentially all the same problem, where for the Stokes problem, as we noted above, the divergence constraint reduces this saddle-point problem to a simpler SPD curl-curl system in terms of its algebraic coarsening, where we assume the block-diagonal saddle-point structure of the Stokes system is given.  

\section{Numerical Results}\label{sec:num}
We present numerical results with Algorithm~\ref{alg:p} applied to curl-curl model problems with different meshes and boundary conditions. We only show the curl-curl results as grad-div cases are similar. We set $\beta = 0.01$, noting that changing $\beta$ typically does not affect convergence. We also demonstrate our solver's effectiveness for the Stokes system. The compatible relaxation we use is the habituated version. We run a two-grid AMG cycle with a pre- and a post-smoothing before restriction and after interpolation, respectively. The convergence tolerance is $10^{-6}$. All problems are solved with a nonzero right-hand side generated randomly.

\begin{figure}[t]
    \centering
    \begin{minipage}{.24\linewidth}    
        \centering
        \begin{tikzpicture}
        \begin{axis}[
            width = 4 cm, height = 6cm,
            font = \tiny,
            ymin = 0, ymax = 1.05,
            xlabel = {number of DoFs},
            scaled x ticks = false,
            legend style={
            font=\tiny,
            at={(0.5,0.815)},
            anchor=north,
            }]
        
            \addplot[color=red, thick, densely dashed]
            table[x=size, y=P, col sep=comma] {plot/convergence_maxwell_periodic_quad.csv};

            \addplot[color=green, mark=star] 
            table[x=size, y=CR, col sep=comma] {plot/convergence_maxwell_periodic_quad.csv};
            
            \addplot[teal, mark =o] 
            table[x=size, y=optimal, col sep=comma] {plot/convergence_maxwell_periodic_quad.csv};

            \addplot[color=magenta, thick] 
            table[x=size, y=P_I0, col sep=comma] {plot/convergence_maxwell_periodic_quad.csv};
            
            \addlegendentry{our $P$}
            \addlegendentry{CR}
            \addlegendentry{$\Popt$}
            \addlegendentry{classical $P_\star$}

        \end{axis}
        \end{tikzpicture}
        \subcaption{uniform quad. }
        \label{fig:plot_I} 
    \end{minipage}
    \hfill
    \begin{minipage}{.24\linewidth}    
        \centering
        \begin{tikzpicture}
        \begin{axis}[
            width = 4 cm, height = 6cm,
            font = \tiny,
            ymin = 0, ymax = 1.05,
            xlabel = {number of DoFs},
            legend style={
            font=\tiny,
            at={(0.5,0.8)},
            anchor=north,
            }]
        
            \addplot[color=red, thick, densely dashed] table[x=size, y=P, col sep=comma] {plot/convergence_maxwell_triangle.csv};
            \addplot[color=blue, mark=triangle]
            table[x=size, y=Hiptmair, col sep=comma] {plot/convergence_maxwell_triangle.csv};
            \addplot[color=green, mark=star] 
            table[x=size, y=CR, col sep=comma] {plot/convergence_maxwell_triangle.csv};
            \addplot[color=magenta, thick] 
            table[x=size, y=P_classical_ideal, col sep=comma] {plot/convergence_maxwell_triangle.csv};

             \addlegendentry{our $P$}
             \addlegendentry{geometric $P$}
             \addlegendentry{CR}
             \addlegendentry{classical $P_\star$}
        \end{axis}
        \end{tikzpicture}
        \subcaption{uniform triangle. }
        \label{fig:plot_II} 
    \end{minipage}
    \hfill
    \begin{minipage}{.24\linewidth}    
        \centering
        \begin{tikzpicture}
        \begin{axis}[
            width = 4cm, height = 5.03cm,
            font = \tiny,
            ymin = 0,
            xlabel = {number of DoFs},
            legend style={
            font=\tiny,
            at={(0.5,1.35)},
            anchor=north}]
        
             \addplot[color=red, thick, densely dashed] table[x=size, y=step_pcg_P, col sep=comma] {plot/convergence_maxwell_triangle.csv};
             \addplot[color=blue,  mark=triangle] table[x=size, y=step_pcg_hiptmair, col sep=comma] 
             {plot/convergence_maxwell_triangle.csv};
             \addplot[color=brown, mark=|] 
             table[x=size, y=step_cg, col sep=comma] {plot/convergence_maxwell_triangle.csv};
             \addplot[color=magenta] 
             table[x=size, y=step_pcg_P_I0, col sep=comma] {plot/convergence_maxwell_triangle.csv};

             \addlegendentry{our $P$}
             \addlegendentry{geometric $P$}
             \addlegendentry{CG}
             \addlegendentry{classical $P_\star$}

        \end{axis}
        \end{tikzpicture}
        \subcaption{pcg.}
        \label{fig:plot_III} 
    \end{minipage}
    \hfill
    \begin{minipage}{.24\linewidth}     
        \centering
        \begin{tikzpicture}
        \begin{axis}[
            width = 4cm, height = 6cm,
            font = \tiny,
            ymin = 0, ymax = 1.05,
            xlabel = {number of DoFs},
            legend style={
            font=\tiny,
            at={(0.5,0.65)},
            anchor=north}]
        
            \addplot[color=red, thick, densely dashed] 
            table[x=size, y=ideal_P, col sep=comma] {plot/stokes_bench_symmetrized.csv};
            \addplot[color=black, thick, dotted] 
            table[x=size, y=geometric_P, col sep=comma] {plot/stokes_bench_symmetrized.csv};
            \addplot[color=green, mark=star] 
            table[x=size, y=CR, col sep=comma] {plot/stokes_bench_symmetrized.csv};

            \addlegendentry{block $P$}
            \addlegendentry{global $P$ }
            \addlegendentry{CR }
        \end{axis}
        \end{tikzpicture}
        \subcaption{Stokes.}
        \label{fig:plot_IV} 
    \end{minipage}

    \caption{curl-curl and Stokes tests}
\end{figure}
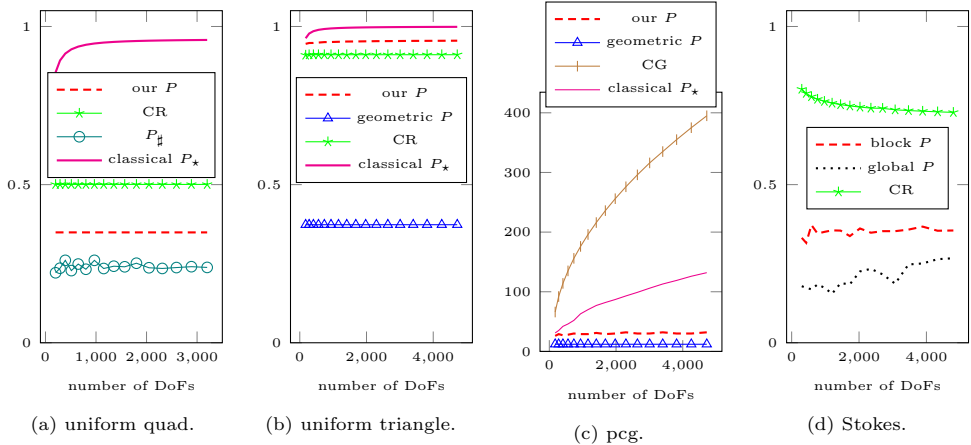

\subsection{Curl-curl experiments}
In Figure~\ref{fig:plot_I}, curl-curl is discretized on a uniform quadrilateral mesh with periodic boundary conditions. For smoothing, we use a global distributive relaxation followed by an L1-Jacobi sweep with a damping factor of 0.5. The smoothers are symmetrized to compare with the optimal interpolation rate. CR is performed on the coarse variables of our algebraic method. For the classical $P_\star$, we run a standard C/F splitting routine on $A$ and form $P_\star$ using $R$ and $S$ as defined in~\eqref{eq:classicalRS}. As expected, the result shows that the convergence rate for classical $P_\star$ diverges to 1. In contrast, our interpolation recovers the N\'{e}d\'{e}lec geometric interpolation on this mesh, and shows a rate only slightly higher than the optimal interpolation $\Popt$~\cite{ali2024generalizedoptimalamgconvergence}, which by theory offers the best possible rate for a given coarsening factor~\cite{ali2024generalizedoptimalamgconvergence}. Additionally, the CR rate demonstrates that the coarse variables selected by Algorithm~\ref{alg:p} form a high-quality set, accurately predicting its good convergence. 

Figure~\ref{fig:plot_II} features curl-curl on the uniform triangular mesh but with Dirichlet boundary conditions. The same smoothers are applied, without symmetrization for simplicity. The convergence rate using the classical $P_\star$ also deteriorates with increasing problem size. While our interpolation by Algorithm~\ref{alg:p} still achieves optimal convergence, the asymptotic rate is slow at 0.94. This is related to our discussion of Figure~\ref{fig:triangle_algebraic} that our algorithm does not fully account for the boundary condition, resulting in missing near-kernel components near the boundary and causing degraded convergence. The high CR rate appropriately predicts the suboptimal choice of the coarse grid and at the same time produces CR error that is large near the boundary. This indicates that we can use global relaxation as well as CR applied globally to $N$ and the error.  Moreover, algebraically, this suggests that a second pass or boundary-specific smoothing could complement our algorithm. However, Figure~\ref{fig:plot_III} demonstrates that our coarse variable construction is on the right track. We apply the same two-grid solver as a preconditioner for the conjugate gradient method (PCG). The plot reaffirms that the classical $P_\star$ is not a good interpolation even with PCG. In comparison, our method is quick to converge with PCG with 30 asymptotic iterations, approximately double that of the geometric interpolation. This indicates that our interpolation captures nearly all of the near-kernel modes that are not sufficiently smoothed by the smoother. Overall, this test validates the correctness and effectiveness of our algebraic construction for preserving the near-kernels on the coarse grid.

\subsection{Stokes experiments}
In Figure~\ref{fig:plot_IV}, the Stokes system is discretized with the MAC scheme on the uniform quadrilateral mesh with periodic boundary conditions. We apply a single sweep of the overlapping Schwarz method (Vanka patch smoother) defined in~\eqref{eqn-g-smoother} for pre- and post-smoothing. First, our two-grid method employs the block interpolation defined in~\eqref{eq:A saddle}, where $P_e$ is the same interpolation our algorithm constructed for the test in Figure~\ref{fig:plot_I}. The uniform convergence rate confirms the close relation between the $\Hcurl$ and Stokes schemes. Next, we form a global interpolation $P$ in the ideal interpolation form by assembling a global $S$ and $R$ that combines the fine and coarse variables defined in the $S$'s and the $R$'s for $P_e$ and $P_n$, respectively. This test further supports the claim of using a stable discretization of the Stokes system using N\'{e}d\'{e}lec elements as discussed in Section~\ref{sec:stokes}. We also consider a sparse approximation to $P$ for the Stokes problem, with results shown in Figure~\ref{fig:stokes_approx} under the same discretization and boundary conditions. We form the block interpolation in~\eqref{eq:A saddle} but instead replaced $P_e$ by its tentative $R^T$. This interpolation directly preserves the constant and can be viewed as a sparse approximation to $P_\star$. We also construct an alternative $R^T$ that forms length-two paths of velocity DoFs connecting adjacent pressure nodal DoFs in the diagonal directions. The corresponding second $P$ is labeled ``(w/ diagonal)'' in Figure~\ref{fig:stokes_approx}. 
\begin{wrapfigure}{r}{0.5\textwidth}
    \centering
    \begin{tikzpicture}
        \begin{axis}[
            width = 4.7cm, height = 4.7cm,
            font = \tiny,
            ymin = 0, ymax = 0.4,
            xlabel = {DoFs},
            ylabel = {convergence rate},
            scaled x ticks = false,
            legend style={
                font=\tiny,
                at={(0.94,0.66)},
                anchor=north}]
        
            \addplot[color=green, mark = star] table[x=Dof,y=no_P_smooth, col sep=comma] {plot/jacob_stokes_bench.csv};
            \addplot[color=blue, mark = star] table[x=Dof,y=P_smooth, col sep=comma] {plot/jacob_stokes_bench.csv};
            \addplot[color=red, thick, densely dashed] table[x=Dof,y=diagonal_no_P_smooth, col sep=comma] {plot/jacob_stokes_bench.csv};
            \addplot[color=brown, mark = triangle] table[x=Dof,y=diagonal_P_smooth, col sep=comma] {plot/jacob_stokes_bench.csv};
             \addlegendentry{unsmoothed $R^T$}
             \addlegendentry{smoothed $R^T$}
             \addlegendentry{unsmoothed $R^T$ (w/ diagonal)}
             \addlegendentry{smoothed $R^T$ (w/ diagonal)}
        \end{axis}
    \end{tikzpicture}
    \caption{Stokes: sparse approximation $P$}
    \label{fig:stokes_approx}
\end{wrapfigure}
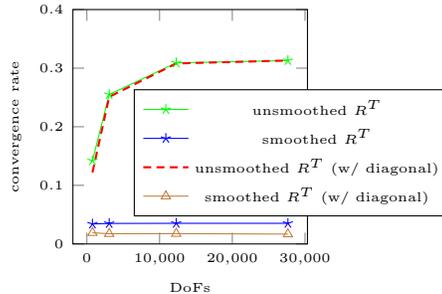Numerical tests show that in this case, gradients are preserved on the coarse grid, and the coarse grid error correction remains divergence-free. The cases where $R^T$'s are smoothed with Richardson are also provided, which results in a solver that resembles the approach from~\cite{voronin2024monolithicalgebraicmultigridpreconditioners}, except that we now coarsen the pressures nodally. 
The results show that unsmoothed $R^T$'s are already effective and that the use of one iteration of Richardson dramatically improves convergence.   The operator complexities for the two-level methods are $1.58, 2.76, 1.79$ and $4.29$, respectively (in legend order).  The high complexities observed for the cases where we match the diagonals are expected because the number of the velocities remains the same on the coarse grid.  From an AMG viewpoint, this suggests that in general using interpolation that preserves gradient structures in Stokes system can benefit the stability on coarse levels.  

\section{Conclusions}\label{sec:conc}
    Overall, we have reduced coarsening the targeted PDE systems to finding a suitable nodal coarsening of the fine-level system in the usual AMG way and then assigning the appropriate averages that preserve the near kernels. Notably, the quality of the coarse grid can depend on the quality of the nodal coarsening, and thus we focus our future efforts on fine-tuning this process. The general algebraic approach that we will explore is using CR to identify inadequately coarsened regions (e.g., near boundaries) and add coarse variables iteratively until CR achieves uniform convergence.  Alternatively, we can explore additional smoothing near the boundaries to avoid coarsening these regions altogether. Finally, we note that the overall approach we propose is a straightforward application of existing ideas in a new way using AMG principles that we hope will provide insights into the AMG coarsening process for PDE systems and other challenging problems.

\end{document}